\documentclass{article}
\usepackage{amsthm,amsfonts,amsmath}
\input epsf.tex
\newdimen\epsfxsize
\newdimen\epsfysize
\def\qed{\vrule height5pt width3pt depth.5pt}

\theoremstyle{plain}
\newtheorem{thm}{Theorem}[section]

\newtheorem{lem}[thm]{Lemma}
\newtheorem{prop}[thm]{Proposition}

\newtheorem{rem}{Remark}[section]

\begin{document}

\title{Virtual Knot Diagrams and the Witten-Reshetikhin-Turaev Invariant}         
\author{H. A. Dye \\
MADN-MATH \\
United States Military Academy \\
646 Swift Road \\
West Point, NY 10996 \\
hdye@ttocs.org \\
Louis H. Kauffman \\
Department of Mathematics, Statistics, and Computer Science  \\
  University of Illinois at Chicago \\
  851 South Morgan St \\
 Chicago, IL  60607-7045 \\ 
 kauffman@uic.edu }        
\date{}          
\maketitle

\begin{abstract}

The Witten-Reshetikhin-Turaev invariant of classical link diagrams is 
generalized to virtual link diagrams. This invariant is unchanged by the framed Reidemeister moves and the Kirby calculus. As a result, it is also an invariant of the 3-manifolds represented by the classical link diagrams. This generalization is used to demonstrate that there are virtual knot diagrams with a non-trivial Witten-Reshetikhin-Turaev invariant and trivial 3-manifold fundamental group. 
\end{abstract}

\section{Introduction} \label{intro}

Generalizations of the Jones polynomial have been defined for virtual link diagrams \cite{kvirt}, \cite{kd1},   \cite{vom}. We continue this process by generalizing the Witten-Reshetihkin-Turaev invariant and the colored Jones polynomials to virtual link diagrams. The Witten-Reshetihkin-Turaev invariant of classical link diagrams is unchanged by the framed Reidemeister moves.  It is also 3-manifold invariant and is unchanged by the Kirby calculus. This invariant is a weighted sum of colored Jones polynomials. Our generalization will also be a weighted sum of colored Jones polynomials and invariant under the framed Reidemeister moves and a generalization of the Kirby calculus.

Every framed, classical link diagram represents a 3-manifold via framed surgery performed on that link. 
Two 3-manifolds represented by link diagrams are homeomorphic if and only if their diagrams are related by a sequence of framed Reidemeister moves and Kirby calculus moves. The fundamental group of this 3-manifold can be computed from a corresponding link diagram. As a result, we can formally extend the definition of the 3-manifold group to a virtual link diagram. We define a virtual Kirby calculus by generalizing the moves of Kirby calculus to incorporate virtual crossings.

The 3-manifold group of a virtual link diagram is invariant under the framed Reidemeister moves, the virtual Reidemeister moves, and the virtual Kirby calculus. This invariance gives rise to the concept of `virtual 3-manifolds'. Virtual link diagrams appear to have the same structure as classical link diagrams (that specify 3-manifolds), but virtual 3-manifolds (specified by non-classical virtual diagrams) have no  corresponding surgery construction. A virtual 3-manifold is, by definition, an equivalence class of link diagrams related by the framed Reidemeister moves, virtual Reidemeister moves, and virtual Kirby calculus.

In this paper, we recall virtual knot theory in section \ref{vkd}. We recall the definitions of the virtual fundamental group and the generalized Jones polynomial. Then, we define the 3-manifold group of a virtual 
link diagram using the definition of the virtual fundamental group. 
We prove that the 3-manifold group of a virtual link diagram is invariant under the framed Reidemeister moves, virtual Reidemeister moves and the Kirby calculus. 

In section \ref{wrti}, we recall the definition of the colored Jones polynomial of a virtual link diagram. We use this definition to define the Witten-Reshetihkin-Turaev invariant of a virtual link diagram. We introduce several formulas that reduce the computational complexity of the Witten-Reshetihkin-Turaev invariant.  Finally, we prove that the Witten-Reshetihkin-Turaev invariant is invariant under the Kirby calculus.

The reader should note that all formulas in the paper that are written in the form of a graphical equation are formulas that indicate the possibility of a substitution of one type of graphic for another graphic in a link diagram. The equality in the equation means that the bracket evaluation of the two diagrams are equal. 

We compute the Witten-Reshetihkin-Turaev invariant of two virtual knot diagrams in the next section. The fundamental group of these virtual knot diagrams is $ \mathbb{Z} $. The 3-manifold group of one of the virtual knot diagrams is the trivial group. However, the Witten-Reshetihkin-Turaev invariant of this virtual knot diagrams is non-trivial. This produces a virtual counterexample to the Poincare conjecture.

In the final section, we discuss the difficulties of computing the Witten-Reshetihkin-Turaev invariant in
 the virtual case. We observe that in order to compute the Witten-Reshetihkin-Turaev invariant of a virtual link diagram, we must compute the Jones polynomial of some complicated virtual link diagrams. The formulas given in sections \ref{wrti} and \ref{wrtexa} were computed using Temperly-Lieb recoupling theory. The full machinery of recoupling theory does not generalize to the virtual case. We present several propositions that highlight the difficulty of computing formulas that include virtual crossings.

\section{Virtual Knot Diagrams} \label{vkd}

Virtual knot theory is a generalization of classical knot theory introduced by Louis H. Kauffman in 1996 \cite{kvirt}. In this section, we review classical and virtual knot theory. We then recall several classical knot invariants that have been generalized to virtual knot diagrams.
 
\begin{figure}[hbt] \epsfysize = 0.65 in
\centerline{\epsffile{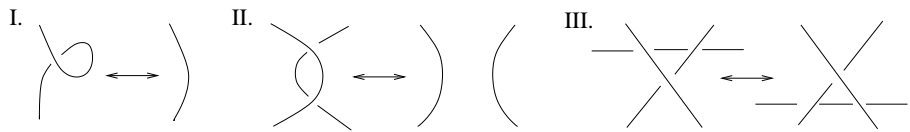}}
\caption{Reidemeister moves}
\label{fig:rmoves}
\end{figure}

A \emph{classical knot diagram} is a decorated immersion of $ S^1 $ into the plane with over/under markings at each crossing. Two classical knot diagrams are said to be \emph{equivalent} is one may be transformed into another by a sequence of Reidemeister moves. Local versions of the Reidemeister move are illustrated in figure \ref{fig:rmoves}.

\begin{figure}[hbt] \epsfysize = 1.5 in
\centerline{\epsffile{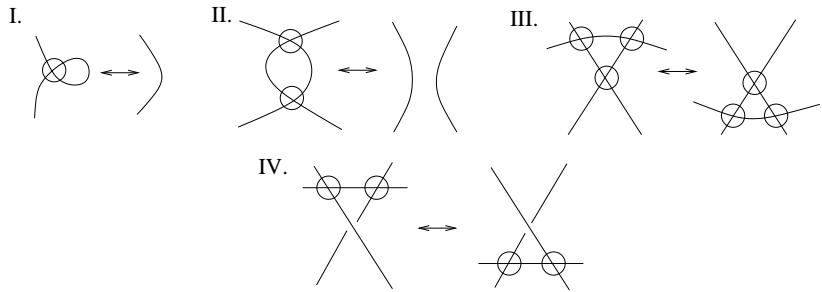}}
\caption{Virtual Reidemeister moves}
\label{fig:vrmoves}
\end{figure}

A \emph{virtual knot diagram} is a decorated immersion of $ S^1 $ into the plane with two types of crossings: classical and virtual. Classical crossings are indicated by over/under markings and virtual crossings are indicated by a solid encircled X. Two virtual knot diagrams are said to be \emph{equivalent} if one may be transformed into the other via a sequence of Reidemeister moves and virtual Reidemeister moves. The \emph{virtual Reidemeister moves} are illustrated in figure \ref{fig:vrmoves}.
Note that the classical knot diagrams are a subset of the virtual knot diagrams.

\begin{figure}[hbt] \epsfysize = 0.8 in
\centerline{\epsffile{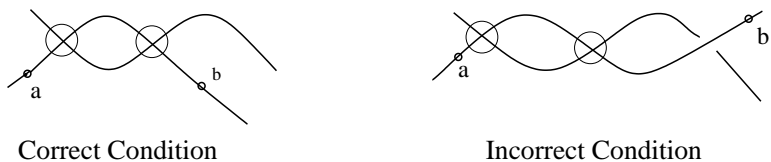}}
\caption{The detour move condition}
\label{fig:detmove}
\end{figure}

The virtual Reidemeister moves can be constructed from a single diagrammatic move: the \emph{detour move}. The detour move pertains only to virtual crossings
and is performed in the following way. Choose an orientation of the virtual knot diagram for
convenience. Select two points, $ a $ and $ b $, on
the virtual knot diagram. The arc from $ a $ to $ b $,
in the direction of the orientation, contains only virtual crossings and no classical crossings. Otherwise we may not perform the detour move.  
(This condition is illustrated in figure \ref{fig:detmove}.) Remove the arc of the knot diagram from point $ a $ to point $ b $, and insert any 
new arc that does not create triple points.
 All double points produced by the placement of the new arc
 are virtual crossings.

\begin{prop}
The detour move is equivalent to a sequence of  virtual Reidemeister moves.
\end{prop}
\begin{rem}We sketch the proof of this proposition. A large scale move like the detour move may be broken down in a series of smaller sequential moves. These smaller moves involve at most three crossings. Note that each virtual Reidemeister move is a detour move.
\end{rem}

\begin{rem}Equivalence classes of virtual knot 
diagrams are in one to one 
correspondence with  equivalence classes of embeddings of 
$ S^1 $ into thickened, oriented 2-dimensional surfaces \cite{kamada}, \cite{kup} and \cite{cks}. A surface with an embedded $ S^1 $ is equivalent to a (possibly different) surface with an embedded $ S^1 $ if one surface with an embedding can be transformed into the other by a sequence of homeomorphisms of the surface, handle cancellations and additions (to the surface), and Reidemeister moves in the surface.
\end{rem}

\begin{figure}[hbt] \epsfysize = 1.0 in
\centerline{\epsffile{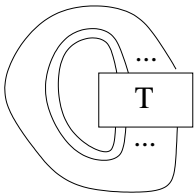}}
\caption{The closure of an n-tangle}
\label{fig:tangleclosure}
\end{figure}

A \emph{ virtual $ n $-tangle diagram} is a decorated immersion of n copies of $ [0,1] $ with classical and virtual crossings. The $ 2n$ endpoints of the diagram are arranged so that $ n$ 
endpoints appear in a row at the top  and $ n $ endpoints form a lower row. Note that it is possible to construct tangle diagrams with no crossings. 
The \emph{closure} of an $ n $-tangle, $ T $ connects each upper endpoint with a corresponding lower endpoint as illustrated in figure \ref{fig:tangleclosure}.

\begin{figure}[hbt] \epsfysize = 0.75 in
\centerline{\epsffile{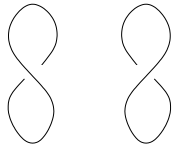}}
\caption{$ \pm 1 $ Framed unknots}
\label{fig:framed}
\end{figure}

We now recall the definition of several invariants of virtual knot diagrams. These invariants will be referred to throughout this paper. 
The \emph{writhe} of a virtual knot diagram and the \emph{linking number} of a virtual link diagram
are computed from an oriented diagram. To define these invariants, we introduce the \emph{sign (crossing sign)} of a classical crossing. 
\begin{figure}[hbt] \epsfysize = 1.0 in
\centerline{\epsffile{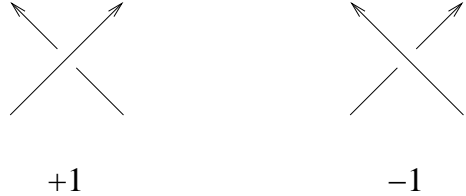}}
\caption{Crossing sign}
\label{fig:sgn}
\end{figure}
Each classical crossing in a virtual knot diagram can be assigned a \emph{crossing sign}. To compute the crossing sign,
we orient the knot diagram. Each classical crossing is assigned a $ \pm $ value based on the orientation 
as shown in  figure \ref{fig:sgn}. Let $ v $ be a classical crossing in a virtual knot diagram. We denote the crossing sign of $ v $ as $ sgn(v) $.
The \emph{writhe} of a virtual knot diagram $ K $ is determined by the crossing sign of each classical crossing in the diagram. 
We denote the writhe of $ K $ as $ w(K) $.
\begin{equation*}
w(K) = \underset{ v \in K}{\sum} sgn(v)
\end{equation*}
The Reidemeister I move changes the writhe of a knot diagram. An unknot with a single Reidemeister I move, is referred to as $ \pm 1 $ framed unknots based on the sign of the crossing.  (We
 illustrate $ \pm 1 $ framed unknots in  figure \ref{fig:framed}.) 
The Reidemeister II and III moves do not change the writhe of a knot diagram and are referred to as the \emph{framed Reidemeister moves}. Two classical knot diagrams that are related by a sequence of Reidemeister II and III moves are said to be \emph{regularly isotopic}. 

The \emph{linking number} of an oriented two component virtual link diagram with oriented components, $ L_1 $ and $ L_2 $, is denoted $ lk (L_1, L_2 ) $. 
Let $ U $ denote the set of all classical crossings formed by the components $ L_1 $ and $ L_2 $. Then the linking number is defined by the formula:
\begin{equation}
 lk(L_1, L_2) = \underset{ v \in U }{\sum} \frac{sgn(v)}{2}. 
\end{equation}
Note that $ lk(L_1, L_2) = lk(L_2, L_1 ) $.
\begin{rem} There is an alternative definition of writhe 
where $ lk( L_1, L_2) $ and $ lk(L_2, L_1 ) $ are not necessarily equivalent \cite{GPV}.
\end{rem}

We recall the fundamental group and the longitude of a virtual knot diagram. These definitions are used to generalize the 3-manifold group obtained from a classical knot diagram to virtual knot diagrams. 

\begin{figure}[hbt] \epsfysize = 1.0 in
\centerline{\epsffile{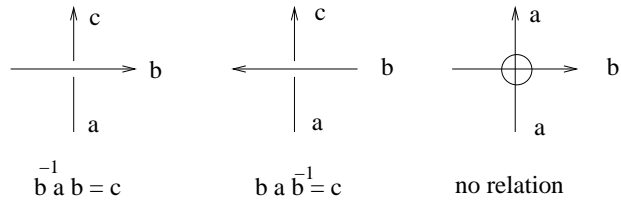}}
\caption{Relations determined by crossings}
\label{fig:fund}
\end{figure}

The \emph{fundamental group} of a virtual knot diagram is a free group modulo relations determined by the classical crossings in the diagram. Each \emph{segment} of the diagram starts and terminates at classical crossings (possibly passing through virtual crossings). We assign each segment a label. These labels represent the generators of the free group. From each crossing we obtain a relation as illustrated in figure \ref{fig:fund}.

The fundamental group of a virtual knot diagram is invariant under the Reidemeister moves and virtual Reidemeister moves.

\begin{figure}[hbt] \epsfysize = 1.0 in
\centerline{\epsffile{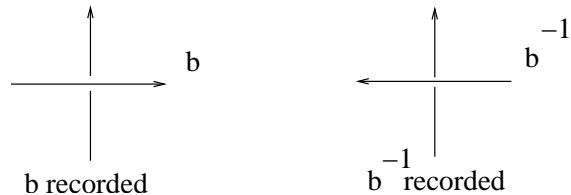}}
\caption{Longitude determined by crossings}
\label{fig:long}
\end{figure}

A \emph{longitude} of a virtual knot diagram is word of the form $ l_1^{ \pm 1} l_2^{ \pm 1} \dots l_1^{ \pm 1} $
where each $ l_i $ is one of the generators of the fundamental group. 
We compute a longitude by choosing an initial point on the  diagram. We traverse the virtual knot diagram from the initial point in the direction of its orientation. We record either $ b $ or $ b^{-1} $ as we traverse an underpass under the arc $ b $. We determine whether to record $ b $ or $ b^{-1} $  as shown in figure \ref{fig:long}. The longitude of a virtual knot diagram is determined up to cyclic permutation of its members.

The \emph{three-manifold group} of a virtual knot diagram $ K $, $ \pi_M (K) $, 
 is the fundamental group modulo the longitude. The three-manifold group is invariant under the Reidemeister moves II and III, the virtual Reidemeister moves, and the Kirby calculus. 
The \emph{virtual Kirby calculus} consists of two moves: the introduction and deletion of $ \pm 1 $ framed unknots and handle sliding.
We describe a \emph{handle slide} for a framed virtual link diagram
$ K $ with at least two components, $ A $ and $ B$. 
Diagrammatically, handle slide $ A $ over $ B $ by the following 
process. Replace $ B $ with two parallel copies of $ B $.
Take the connected sum of $ A $ and the outermost copy of $ B $.

\begin{figure}[hbt] \epsfysize = 1.0 in
\centerline{\epsffile{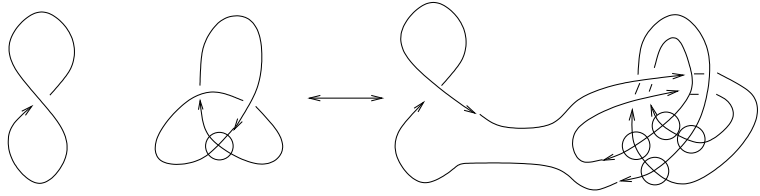}}
\caption{A example of a handle slide with virtual crossings}
\label{fig:slideexa}
\end{figure}

We show an example of this process in  figure \ref{fig:slideexa},
where a $ + 1 $ framed unknot is slid over a trefoil.

These two moves are called the \emph{Kirby moves}.
\begin{rem} In the case of classical link diagrams, the definition of the classical Kirby calculus is identical to the definition given above without virtual crossings. 
\end{rem}
\begin{rem} The following question has been posed: if two classical diagrams are equivalent under the virtual Kirby calculus, are they equivalent under the classical Kirby calculus.
\end{rem}
The 3-manifold group is invariant under the framed Reidemeister moves and the virtual Reidemeister moves. It remains to prove that the 3-manifold group is invariant under the Kirby moves.

A classical framed link diagram represents the three manifold obtained by 
removing a small neighborhood of each component of the link diagram, and gluing in 
a solid torus so that the meridian is identified with the
longitude of the removed link component \cite{purp}. 
The Reidemeister I move changes the longitude of a knot diagram and changes the identification used to glue in the torus.

For example, the fundamental group of an unknot with writhe $ 0 $ and the fundamental group of a  $ \pm 1 $ framed unknot are both $ \mathbb{Z} $. However, the 3-manifold group of an unknot with writhe zero is $ \mathbb{Z} $ and  the three-manifold group of 
the $ \pm 1 $ framed unknot is the trivial group. The unknot with writhe $ 0 $ represents $ S^1 \times S^2 $ and 
the $ \pm 1 $ framed unknot represents $ S^3 $.

\begin{thm}
If two classical framed link diagrams, $ K $ and $ \hat{K} $ , are related 
by a sequence of Kirby moves and framed Reidemeister moves then 
$$
  \pi_M (K) = \pi_M ( \hat{K} ) 
$$
and the two diagrams represent homeomorphic three manifolds.
\end{thm}
\textbf{Proof: } \cite{purp}. \qed

\begin{thm} The 3-manifold group of a virtual knot diagram $ K $, 
 is invariant
under the framed Reidemeister moves, the virtual Reidemeister moves (detour move) and the Kirby
calculus.
\end{thm}

\textbf{Proof:}
 The fundamental group is not changed by the framed Reidemeister moves and
 the virtual Reidemeister moves. The longitude is also unchanged by 
 the Reidemeister II and III moves and the detour move. This implies 
 that the 3-manifold group of a virtual knot diagram is not changed by
 these moves. 
 We only need to show that the 3-manifold group is unchanged by 
 the Kirby moves. 

 The 3-manifold group of a $ \pm 1 $  framed unknot is 
 trivial and
 the  introduction of a  disjoint $ \pm 1 $  framed unknot
 to the diagram does not change the 3-manifold group.

\begin{figure} [hbt] \epsfysize = 0.5 in
\centerline{\epsffile{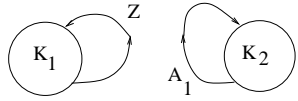}}
\caption{Two virtual knot diagrams before a handle slide}
\label{fig:init}
\end{figure}

We consider the effect of handle sliding on a virtual knot diagram.
Consider a virtual knot diagram with two disjoint 
components, oriented as shown in  figure \ref{fig:init}.   Denote the two
components by
 $ K_1 $ and $ K_2 $, the 3-manifold group by 
$ \pi_M (K_1  \amalg K_2 ) $, and  the
longitudes of $ K_1 $ and $ K_2 $ as $ L_1 $ and $ L_2 $.   
Assume that the arcs of $ K_1 $ are labeled sequentially 
 as $ Z, Z_1, Z_2 \ldots Z_m $ in the direction of orientation.
The placement of the labeled arc $ Z $ is indicated in  figure \ref{fig:init}.
 Similarly,  the arcs of the component $ K_2 $ are labeled in
the direction of orientation as $ A_1 ,  A_2 , \ldots A_n $.
(To generalize the argument, we would
consider the different possible orientations of $ K_1 $ and $ K_2 $, 
 and the case where $ K_2 $ 
is a link with multiple components.)

\begin{figure} [hbt] \epsfysize = 1.0 in
\centerline{\epsffile{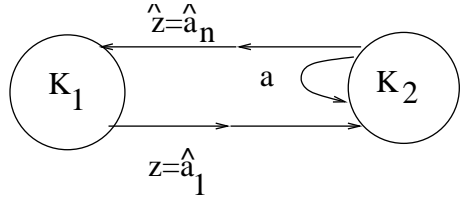}}
\caption{Two virtual knot diagrams after a handle slide}
\label{fig:hstate2}
\end{figure}

Handle slide $ K_1 $ over $ K_2 $, and we obtain the diagram 
shown in  figure \ref{fig:hstate2}. Let $ \hat{K} $ denote the new diagram,
$ \pi_M (\hat{K}) $ the corresponding 3-manifold group, and $ \hat{L_1} $
and $ \hat{L_2 } $ the new longitudes.

\begin{figure}[hbt] \epsfysize = 1.0 in
\centerline{ \epsffile{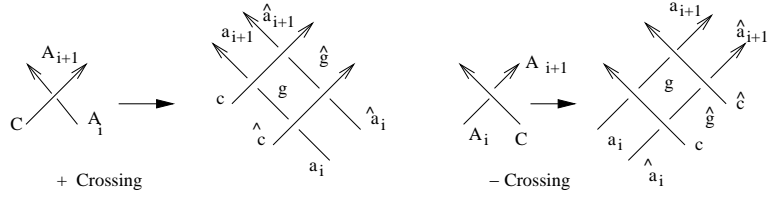}}
\caption{Detail of of a crossing before and after a handle slide}
\label{fig:hstate3}
\end{figure}

The virtual crossings do not produce any new relations, therefore 
we need only examine the effect of the handle slide on the classical crossings. 
Each original crossing in $ K_1 $ is unchanged. 
Each classical crossing in $ K_2 $ becomes a set of four classical crossings whose four relations reduce to
a single equation. 
Consider a positive crossing in $ K_2 $ as shown in figure
\ref{fig:hstate3}.

We obtain 
the following relationship from a positive classical
crossing before handle sliding:
$
 C^{-1} A_i C = A_{i+1}$.
After handle sliding, we obtain four relations from the positive crossing.
\begin{alignat*}{2}
c^{-1} g c &= a_{i+1} & \qquad
c^{-1} \hat{g} c &= \hat{a}_{i+1}\\
\hat{c}^{-1} a_i \hat{c} &= g & \qquad
\hat{c}^{-1}\hat{a}_i \hat{c} &= \hat{g} 
\end{alignat*}
These equations may be simplified to two equations:
\begin{align*}
c^{-1}\hat{c}^{-1} a_i \hat{c}   c &= a_{i+1} \\
c^{-1}\hat{c}^{-1}\hat{ a}_i \hat{c} c &= \hat{a}_{i+1}
\end{align*} 
Multiplying, we obtain a single relation from each set of
four  positive crossings:
$$
(\hat{c} c)^{-1} ( \hat{a}_i a_i)( \hat{c} c) = \hat{a}_{i+1} a_{i+1}
$$
The change for negative classical crossing is also shown in   figure
 \ref{fig:hstate3}. We also obtain a single relation in the case of 
a negative crossing.
$$
(\hat{c} c)( \hat{a}_i a_i)  ( \hat{c} c)^{-1} =  \hat{a}_{i+1}a_{i+1}
$$

Consider the longitude of $ K_2 $. Compute the original longitude by
traversing the knot starting at $ A _1$. If
 $ L_2 = \lbrace x_1^{\pm 1}, x_2^{\pm 1} \ldots x_n^{\pm 1} \rbrace $
then  
\begin{equation*} \hat{L_2} =
 \lbrace (\hat{x_1}x_1)^{\pm 1} , (\hat{x_2}x_2)^{\pm 1}  \ldots
 (\hat{x_n}x_n)^{\pm 1} \rbrace .
 \end{equation*}
Hence, $ \hat{L_2} a_1 z  \hat{L_2}^{-1} = a_1 \hat{z}$ 
implies that
$ z = \hat{z} $. 
Further, if $ L_1 $ is computed by traversing the knot from z in the 
direction of the orientation, then we observe that 
 $  \hat{L_1} = L_1 \hat{L_2} $. 
These computations allow us to conclude that 
$ \pi_M (K) = \pi_M ( \hat{K} ) $.  \qed

\begin{figure}[htb] \epsfysize = 1 in
\centerline{\epsffile{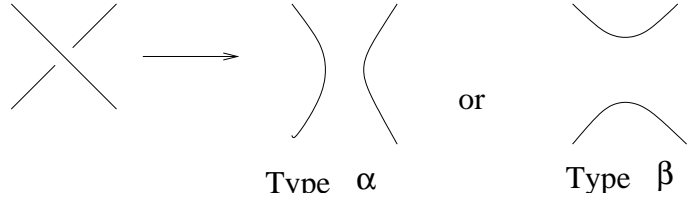}}
\caption{Smoothing types}
\label{fig:smooth}
\end{figure}

We recall a generalization of the bracket model \cite{Kdetect} to define the \emph{generalized Jones polynomial} of a virtual knot diagram.
Each  classical crossing may be \emph{smoothed} as a type 
$ \alpha $ smoothing or  a type $ \beta $ smoothing
as illustrated in figure \ref{fig:smooth}.

A \emph{state} of a virtual knot diagram is a diagram determined by 
a choice of smoothing type for each classical crossing. Note that a state of
a virtual link diagram may contain virtual crossings and consists of closed 
curves (possibly) with virtual crossings.

Let $ K $ be a virtual knot diagram. 
Let $ s $ be a state of the virtual knot diagram, let $ S $ represent the set of 
all possible states and let 
$ c(s) $ equals the number of type $ \beta $ crossings minus the number of type $ \alpha $ crossings. 
Let $ |s| $ denote the number of closed curves in the smoothed diagram and $ d = -(A^2 + A^{2}) $. The bracket polynomial of $ K $ 
denoted $ \langle K \rangle  $ then:
\begin{equation*}
 \langle K \rangle = \underset{s \in S}{\sum} d^{|s|-1} A^{c(s)}
\end{equation*}
The bracket polynomial is invariant under the Reidemeister II and III moves and the virtual Reidemeister moves. To obtain invariance under the Reidemeister I move, we use the writhe of $ K $ to normalize. 
We define the \emph{normalized bracket polynomial} of a virtual knot diagram.
\begin{equation*}
  f_K (A) = (-A)^{-3w(k)} \langle K \rangle 
\end{equation*}
The only difference between this definition and the definition of the classical Jones polynomial is that a state of virtual knot diagram is set of closed curves (possibly containing virtual crossings) and not a set of simple closed curves.
The \emph{Jones polynomial}, $ V_K (t) $,and the bracket polynomial are related by letting $ A = t^{- 1/4}$.
\begin{equation}
 f_K (t^{1/4})= V_K (t)
\end{equation}
In the next section, we introduce the Witten-Reshetihkin-Turaev invariant. This invariant will be generalized to virtual knot diagrams and is invariant under the Kirby calculus, framed Reidemeister moves, and the virtual Reidemeister moves.

\section{The Witten-Reshetikhin-Turaev Invariant} \label{wrti}

In this section, we extend the definition of the Witten-Reshetikhin-Turaev invariant \cite{witten}, \cite{rt1}, \cite{rt2} to virtual link diagrams.   First, we recall the 
definition of the Jones-Wenzl projector (q-symmetrizer)\cite{tl}.
We then define the colored Jones polynomial of a virtual link diagram. 
It will be clear from this definition that two equivalent virtual knot diagrams have the same colored Jones polynomial. We will use these definitions to extend the Witten-Reshetikhin-Turaev invariant to virtual link diagrams. From this construction, we conclude that two virtual link diagrams, related by a sequence of framed Reidemeister moves and virtual Reidemeister moves, have the same Witten-Reshetikhin-Turaev invariant. Finally, we will prove that the generalized Witten-Reshetikhin-Turaev invariant is unchanged by the virtual Kirby calculus. In the next section, we present two virtual knot diagrams that have fundamental group $ \mathbb{Z} $ and a Witten- Reshetikhin-Turaev invariant that is not equivalent to $ 1 $.

To  form the \emph{n-cabling} of a virtual knot diagram, take $ n $ parallel copies of the virtual knot diagram. A single classical crossing becomes a pattern of 
$ n^2 $ classical crossings and a single virtual crossing becomes $ n^2 $ virtual crossings.

Let $ r $ be a fixed integer such that $ r \geq 2 $ and let 
\begin{equation*}
 A = e^{\frac{\pi i}{ 2r } }.
\end{equation*} 
Here is  a formula used in the construction of the Jones-Wenzl projector.  
\begin{equation*} \Delta_n = (-1)^n \frac{A^{2n + 2} - A^{-(2n+2)}}{A^2 - A^{-2}}. 
\end{equation*}
 Note that $ \Delta_1 = - (A^2 + A^{-2})$, the value assigned to a simple closed curve by the bracket polynomial.
 There will be an an analogous interpretation of $ \Delta_n $ which we will discuss later in this section.  
 
\begin{figure}[htb] \epsfysize = 1 in
\centerline{\epsffile{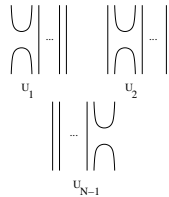}}
\caption{Elementary n-tangles}
\label{fig:tlgen}
\end{figure}

 We recall the definition of an n-tangle. Any two n-tangles can be multiplied by attaching the 
 bottom n strands of one n-tangle to the upper n strands of another n-tangle. We define an n-tangle to be \emph{elementary} if it contains to no classical or virtual crossings. Note that the product of any two elementary tangles is elementary.
 Let $ I $ denote the identity n-tangle and let $ U_i $ such that $ i \in \lbrace 1,2, \dots n-1 
\rbrace $ denote the n-tangles shown in figure \ref{fig:tlgen}. 
By multiplying a finite set of $ U_{i_1} U_{i_2} \dots U_{i_n} $, we can obtain any elementary n-tangle.
Formal sums of the elementary tangles over $ \mathbb{Z} [A, A^{-1} ] $ generate the $ n^{th} $ Temperly-Lieb algebra \cite{tl}.
 
We recall that the \emph{$ n^{th} $ Jones-Wenzl projector} is a certain sum of all elementary $ n $-tangles with coefficients in $ \mathbb{C}$ \cite{tl}, \cite{knotphys}. We denote the $ n^{th} $ Jones Wenzl projector as 
$ T_n $. We indicate the presence of the Jones-Wenzl projector and the n-cabling by labeling the component of the knot diagram with $n$.
\begin{rem}There are different methods of indicating the presence of a Jones-Wenzl projector. In a virtual knot diagram, the presence of the $ n^{th} $ Jones-Wenzl projector is indicated by a box  with n  strands entering and n strands leaving the box. For n-cabled components of a virtual link diagram with a attached Jones-Wenzl projector, we indicate the cabling by labeling the component with $ n $ and the presence of the projector with a box. This notation can be simplified to the convention indicated in the definition of the the colored Jones polynomial. The choice of notation is dependent on the context.
\end{rem}

\begin{figure}[htb] \epsfysize = 1 in
\centerline{\epsffile{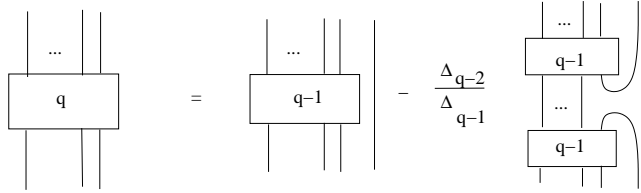}}
\caption{$ n^{th} $ Jones-Wenzl Projector}
\label{fig:asym}
\end{figure}

We construct the Jones-Wenzl projector recursively.  
The $ 1^{st} $ Jones-Wenzl projector consists of a single strand with coefficient $ 1 $. The is exactly one 1-tangle with no classical or virtual crossings. 
The $ n^{th} $ Jones-Wenzl projector is  constructed from the $ (n-1)^{th} $ and $ (n-2)^{th} $ Jones-Wenzl projectors as illustrated in figure \ref{fig:asym}. 
 
\begin{figure}[htb] \epsfysize = 1 in
\centerline{\epsffile{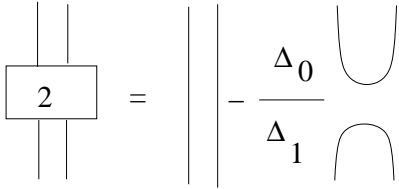}}
\caption{$ 2^{nd} $ Jones-Wenzl Projector}
\label{fig:2sym}
\end{figure}

We use this recursion to construct the $ 2^{nd} $ Jones-Wenzl projector as shown in figure \ref{fig:2sym}.

We will refer to the Jones-Wenzl projector as the J-W projector for the remainder of this paper.

We review the properties of the J-W projector.
Recall that $ T_n $ denotes the $ n^{th} $ J-W projector then
\begin{gather*}
\text{i) }  T_n T_m = T_n  \text{ for $ n \geq m $ } \\
\text{ii) } T_n U_i = 0  \text{ for all $ i $ } \\
\text{iii) } \text{The bracket evaluation of the closure of } T_n = \Delta_n 
\end{gather*}

\begin{rem}The combinatorial definition of the J-W projector  is given in \cite{tl}, p. 15.
Note that \cite{tl} provides a full discussion of all formulas given above. 
\end{rem} 

Let $ K $ be a virtual link diagram with components $ K_1, K_2 \dots K_n $. Fix an integer $ r \geq 2$ and let $ a_1, a_2 \dots a_n  \in  \lbrace 0, 1,2, \ldots r-2 \rbrace $. Let $ \bar{a} $ represent the vector $ ( a_1, a_2, \dots a_n ) $. Fix $ A= e^{\frac{\pi i}{2r}} $ and  $ d= -A^2 - A^{-2} $.
We denote the \emph{generalized $ \bar{a} $ colored Jones polynomial} of $ K $ as $ \langle K^{\bar{a}} \rangle $. 
To compute $ \langle K^{ \bar{a} } \rangle $, we
 cable the component $ K_i $ with $ a_i $ strands  and attach the $ a_i^{th} $ J-W projector  to  cabled component $ K_i $. We apply the Jones polynomial to the cabled diagram with attached J-W projectors. 

The colored Jones polynomial is invariant under the framed Reidemeister moves and the virtual Reidemeister moves. This result is immediate, since the Jones polynomial is invariant under the framed Reidemeister moves and the virtual Reidemeister moves.

\begin{rem}The $ a $-colored Jones polynomial of the unknot is $ \Delta_{a} $. In other words, the 
 Jones polynomial of the closure of the $ a^{th} $ J-W projector is $ \Delta_a $.
\end{rem}

The \emph{generalized Witten-Reshetikhin-Turaev invariant} of a virtual link diagram is a sum of 
colored Jones polynomials. 
Let $ K $ be a virtual knot diagram with $ n $ components. Fix an integer $ r \geq 2 $. We denote the unnormalized Witten-Reshetikhin Turaev invariant of $ K $ as $ \langle  K^{ \omega}  \rangle $, which is shorthand for the following equation. 
\begin{equation}
\langle K^{\omega} \rangle = \underset{\bar{a} \in \lbrace 0,1,2, \dots r-2 \rbrace^n }{\sum}
\Delta_{a_1} \Delta_{a_2} \dots \Delta_{a_n} \langle K^{ \bar{a}} \rangle
\end{equation}

\begin{rem}For the remainder of this paper, the Witten-Reshetikhin-Turaev invariant will be referred to as the WRT.
\end{rem}

We define the matrix $ N $ in order to construct the normalized WRT \cite{tl}. 
Let $ N $ be the matrix defined as follows:
\begin{align*}
\text{i) } N_{ij} &=  lk (K_i,K_j) \text{ for } i \neq j \\
     \text{ii) }N_{ii} &= w(K_i) \\
\text{then let} \\
        b_{+} (K) &= \text{ the number of positive eigenvalues of } N \\
        b_{-} (K) &= \text{ the number of negative eigenvalues of } N \\ 
        n(k) &=  b_{+} (K) - b_{-} (K). 
\end{align*}
The normalized WRT of a virtual link diagram $ K $ 
is denoted as $ Z_{K} (r) $.
Let $ A = e^{ \frac{ \pi i }{ 2r}} $ and let $ |k|  $ denote the number of
components in the virtual link diagram $ K $. Then $ Z_K (r) $ is defined 
by the formula
\begin{equation*}
        Z_{K}(r) = \langle K^{ \omega}  \rangle \mu^{ | K | + 1}
        \alpha^{-n(K)}
\end{equation*}
where
\begin{align*}
  \mu &= \sqrt{ \frac{2}{r}} sin (\frac{\pi}{ r}) \\
\text{ and } \\
         \alpha &= (-i)^{r-2} e^{i \pi [ \frac{3(r-2)}{4r}]}.
\end{align*}
This normalization is chosen so that normalized WRT of the  
unknot with writhe zero is $ 1 $ and the normalization is invariant under the introduction and deletion of $ \pm 1 $ framed unknots.

Let $ \hat{U} $ be a $+1 $ framed unknot. We recall that $ \alpha = \mu \langle \hat{U}^{\omega} \rangle $ \cite{tl}, page 146. Since $ \hat{U} $ and $ K $ are disjoint in $ K \amalg \hat{U} $ then
 $ \langle (K \amalg \hat{U} )^{ \omega} \rangle = \langle K^{\omega} \rangle \langle \hat{U}^{\omega} \rangle $. We note that $ b_+ ( K \amalg \hat{U}) = b_+ (K) + 1 $,
 $ b_- ( K \amalg \hat{U}) = b_- (K) $, and $ | K \amalg \hat{U} | = |K|+1 $. We compute that
\begin{equation*}
Z_{K \amalg \hat{U} } (r) = \langle K^{ \omega} \rangle \langle \hat{U}^{\omega} \rangle \mu^{ |K|+2 } \alpha^{-n(K)-1}
\end{equation*}
As a result, 
\begin{equation*} Z_{K \amalg \hat{U} } (r) = Z_K (r).
\end{equation*}

We demonstrate that the normalized, generalized WRT is invariant under the framed 
Reidemeister moves and the Kirby calculus. The WRT is a sum of colored Jones polynomials and it is  clear that the WRT is invariant under  the framed Reidemeister moves and the virtual Reidemeister moves. In particular, the normalized WRT is invariant under the first Kirby move (the introduction and deletion of $ \pm 1 $ framed unknots) due to the choice of normalization. We only need to show invariance under handle sliding.

\begin{thm}\label{invar} Let $ K $ be a virtual link diagram then $ Z_K(r) $ is invariant under the framed Reidemeister moves, virtual Reidemeister moves, and the virtual Kirby calculus. 
\end{thm}
In the next few subsections, we introduce some necessary machinery.

In the classical case, the computation of this invariant is simplified by the use of recoupling theory. In the next section, we introduce formulas from recoupling theory \cite{tl}. Recoupling theory establishes a relationship between labeled knot diagrams and labeled trivalent graphs with classical and virtual crossings.
These formulas allow us to compute the WRT without a full expansion of the Jones polynomial.  

\subsection{Kirby Calculus Invariance} 

 These formulas are obtained from recoupling theory at roots of unity \cite{tl}. The formulas 
simplify the computation of the WRT. They will be used to prove that the WRT is invariant under Kirby calculus.

Consider a virtual link diagram with components: $ K_1, K_2 \dots K_n $. We cable component $ K_i $ with $ a_i $ strand and attach a J-W projector. This virtual link diagram can be regarded as a collection of labeled trivalent graphs  (each with one edge and no vertices) with classical and virtual crossings. We introduce some formulas obtained from recoupling theory \cite{knotphys}. These formulas describe identities between sums of bracket evaluations of labeled trivalent graphs.
A full derivation of each formula can be obtained in \cite{tl}. 

\begin{figure}[htb] \epsfysize = 1.5 in
\centerline{\epsffile{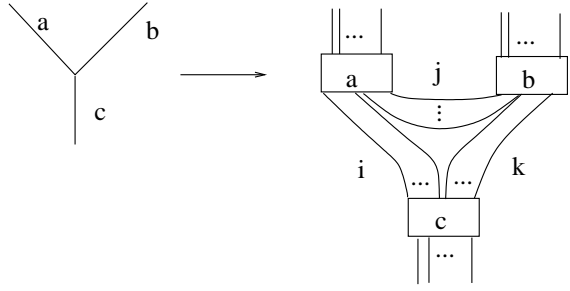}}
\caption{Associating projectors to a trivalent vertex}
\label{fig:assocline}
\end{figure}

Given a trivalent vertex with labels $ a, b$, and $ c $ coming into the vertex, we translate this into an expression of J-W projectors. We connect the internal lines as indicated in figure
\ref{fig:assocline}. Note that $ i + j= a $, $ j + k =b $, and $ k+i=c $ where $ i, j,$ and $ k $ are non-negative integers. In order for such a construction to be made: $ a+ b +c $ is even, $ a + b -c \geq 0 $,  $ b+ c-a \geq 0$, and $ c+ a -b \geq 0 $. These conditions are an admissibilty condition for translating any trivalent vertex into projectors for evaluating the colored Jones polynomial for any value of $ \bar{a} $. We now give a more technical definition of admissibility for evaluation at roots of unity.

Fix an integer $ r \geq 2 $ and let  $ A = e^{\frac{\pi i}{2r}}$.
The \emph{quantum factorial}  is a factorial product based on $ \Delta_{n}$.  Let 
 $ 0 \leq m \leq r-1 $. The $ m^{th} $ quantum 
factorial  is denoted by $ [m]! $. Now,
\begin{gather*}
 [n]  = (-1)^{n-1} \Delta_{n-1}  \\
 [0]! = 1 \text{ and }   [1]! = \Delta_0 \\ 
[m]! =  [m] [m-1]  \ldots [1]. 
\end{gather*} 

For $ r $, an \emph{admissible triple} is a set of three non-negative integers such that:
\begin{gather*}
\text{i)} a + b+ c \leq 2r-4  \text{ and is even} \\
\text{ii) } a + b -c \geq 0 \\
 \text{iii) } a+c -b \geq 0 \\
\text{iv) } b+c -a \geq 0
\end{gather*}

\begin{figure}[hbt] \epsfysize = 1.0 in
\centerline{\epsffile{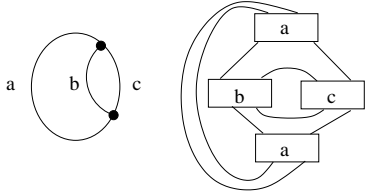}}
\caption{$ \theta $-Net Graph}
\label{fig:tnet}
\end{figure}

We recall the $ \theta $-net formula,  \emph{ $ \theta(a,b,c) $ }. Let $ a,b,c $ be an admissible triple and let
 $ m = \frac{a+b-c}{2} $,
 $ n= \frac{a+c-b}{2}$, and $ p = \frac{b+c-a}{2} $ then 
\begin{equation} \label{thnet} \theta (a,b,c) =
(-1)^{m+n+p} \frac{[m+n+p+1]! [m]! [n]! [p]!}{ [m+n]![n+p]![p+m]!}
\end{equation}
The labeled trivalent diagram shown in figure \ref{fig:tnet} is referred to as a $ \theta $-net.  The virtual link diagram represented by the $ \theta $-net diagram is pictured on the right (by using the prescription for replacing trivalent vertices with projectors and reducing them). The Jones polynomial of this virtual link diagram is $ \theta (a,b,c) $. 

Note that we can only construct the link diagram corresponding to a $ \theta$-net if $ (a,b,c) $ is an admissible triple.

\begin{figure}[hbt] \epsfysize = 1.0 in
\centerline{\epsffile{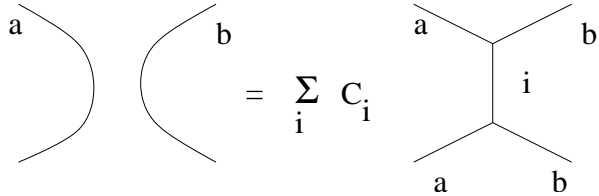}}
\caption{Graph fusion}
\label{fig:fusion}
\end{figure}

The formula shown in figure \ref{fig:fusion} (a special case of the recoupling formulas from section \ref{recoup}) expresses the fact that the bracket evaluation of a diagram with two nearby arcs can be replaced with the bracket evaluation of the sum of diagrams obtained by the indicated graphical insertion. 

The coefficients $ C_i $ in this sum are: 
\begin{equation} \label{fuse}
C_i = \frac{ \Delta_i}{ \theta (a,b,i)}.
\end{equation}
We can only construct the corresponding virtual link diagram when $ (a,b,i) $ is an admissible triple.

These formulas simplify the computation of the WRT of a virtual link diagram.  

To prove that the normalized WRT of a virtual knot diagrams is unchanged by the virtual Kirby calculus, we need only demonstrate that the WRT is unchanged by  handle sliding.
We prove the following lemma:
\begin{lem} \label{hslem} Suppose that $ \langle K^{\omega} \rangle $ satisfies the conditions shown in figure \ref{fig:cond}.
\begin{figure}[htbp]
	\centerline{\epsffile{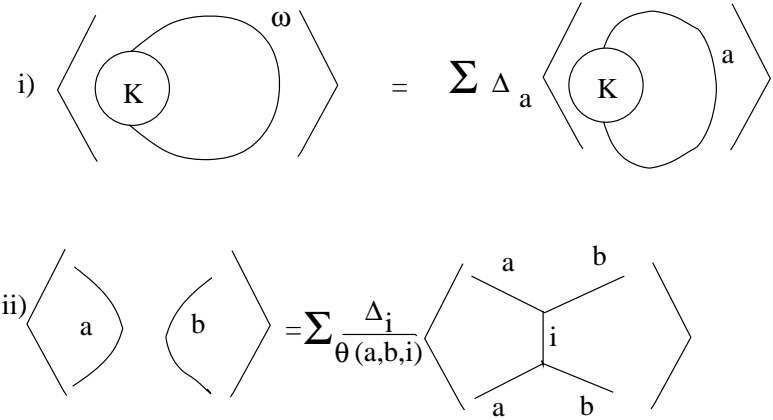}}
	\caption{Conditions for Lemma \ref{hslem}}
	\label{fig:cond}
\end{figure}   
Then $ \langle K^{ \omega } \rangle $ is invariant under handle sliding.
\end{lem}
\textbf{Proof:}

\begin{figure}[hbt] \epsfysize = 1.0 in
\centerline{\epsffile{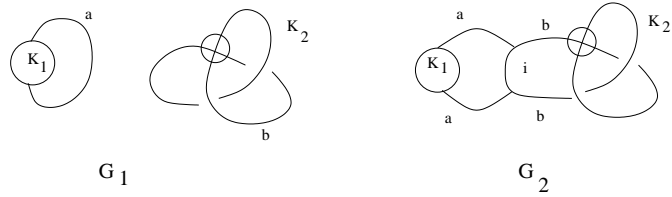}}
\caption{Recoupling move for the handle slide}
\label{fig:g1}
\end{figure}

Let $ K_1 $ and $ K_2 $ be two disjoint components of a virtual knot diagram.
Let $ G_1 $ denote the cabled  $ K \amalg \hat{K}$ with J-W projectors attached  as illustrated in figure \ref{fig:g1}.

We handle slide $ K_1 $ over the knot $ K_2 $ to obtain the diagram $ G_4 $. We will show 
that $ \langle G_1^{ \omega }  \rangle = \langle G_4^{ \omega } \rangle $. 

Fix $ r \geq 2 $ and let
$ I= \lbrace 0,1, \ldots r-2 \rbrace $. Then: 
$$ \label{s1}
 \langle G_1^{\omega}  \rangle = \underset{a, b \in I}{ \sum } 
  \Delta_a \Delta_b \langle ( K_1 \amalg  K_2 )^{( a , b )} \rangle . 
$$

Fix $ a $ and consider the sum:
\begin{equation} \label{s2}
  \underset{b \in I}{ \sum} \Delta_a \Delta_b \langle ( K_1 \amalg K_2 )^{(a,b)}  \rangle 
\end{equation}

Perform fusion on $ K_1 $ and $ K_2 $ and obtain the new diagram $ G_2 $.  
Now:
$$ 
\underset{b \in I}{\sum} \Delta_a \Delta_b \langle G_1^{ (a,b) } \rangle = 
\underset{b \in I}{\sum} \Delta_a \Delta_b \underset{i \in I}{\sum } \frac{ \Delta_i}{ \theta (a,b,i) } \langle G_2^{ (a,b,i) } \rangle
$$

\begin{figure}[hbt] \epsfysize = 1.0 in
\centerline{\epsffile{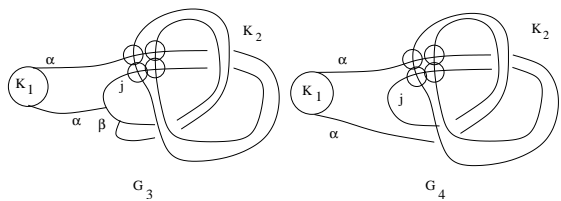}}
\caption{Recoupling diagrams after fusion}
\label{fig:g3}
\end{figure}

We slide the upper vertex of this graph until we obtain the diagram $ G_3 $ as
shown in figure \ref{fig:g3}. Sliding this edge corresponds to a sequence of Reidemeister II, III and virtual Reidemeister moves. 

As a result,
\begin{equation}
\underset{ b \in I}{ \sum } \Delta_a \Delta_b \langle G_1^{ (a,b) }   \rangle = \underset{ b,i \in I}{\sum} \frac{ \Delta_a \Delta_b \Delta_i }{ \theta (a,b,i)} \langle G_3^{ (a,b,i) }  \rangle
\end{equation}
Applying the inverse of the fusion formula to the previous equation: 
\begin{align*}
\underset{ b \in I}{ \sum } \Delta_a \Delta_b \langle G_1^{ (a,b) }   \rangle =
\underset{ i \in I }{ \sum } \Delta_a \Delta_i \langle G_4^{ (a, i) } \rangle
\end{align*}
Now:
\begin{equation}
\underset{ a , b \in I }{ \sum } \Delta_a \Delta_b \langle G_1^{ (a,b) } \rangle  = \underset{a, i \in I}{ \sum} \langle G_4^{ ( a , i ) } \rangle
\end{equation}
Hence, 
$
\langle G_1^{ ( \omega , \omega ) } \rangle = \langle G_4^{ ( \omega , \omega ) } \rangle.
$
\qed

\section{Examples} \label{wrtexa}
\begin{figure}[hbt] \epsfysize = 1.0 in
\centerline{\epsffile{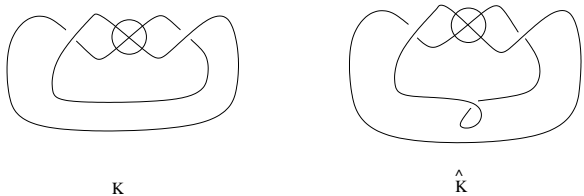}}
\caption{Virtual knot diagrams, $K $ and $ \hat{K}$}
\label{fig:wrtex}
\end{figure}

In this section, we compute the WRT of the two virtual knot diagrams, $ K $ and $ \hat{K} $, shown in figure \ref{fig:wrtex}.

Both $ K $ and $ \hat{K} $ have fundamental group $ \mathbb{Z} $. 
The three-manifold groups of these diagrams are
\begin{equation*}
 \pi_M (K) = \mathbb{Z}_2 \text{ and } \pi_M ( \hat{K} ) = 1 .
 \end{equation*}
Note that $ \pi_M ( S^1 \times S^2 ) = \mathbb{Z} $ and $ \pi_M ( S^3 ) = 1 $.

To compute the WRT of these virtual knot diagrams, we recall the following formulas from recoupling 
theory. These formulas establish a correspondence between labeled trivalent graphs and sums of virtual knot diagrams with coefficients in $ \mathbb{C} $. 

We define the \emph{tetrahedral net} formula. Let $a,b,c,d,e$ and $f $ be non-negative integers.
\begin{equation} \label{tetform}
Tet \begin{bmatrix}  a & b & e \\ c & d & f \end{bmatrix}  = 
 \frac{ G  }{ E } \prod_{m \leq s \leq M}
 \frac{ (-1)^s [s+1]!}{ \prod_{i} [s-a_i]! \pi_j [b_j -s]!} 
\end{equation} 
such that
\begin{alignat*}{2} 
        E &= [a]![b]![c]![d]![e]![f]!  & \qquad
        G &= \prod_{i,j} [b_j - a_i]! \\
        a_1 &= \frac{( a +d +e) }{2}  & \qquad
        a_2  &= \frac{(b+c+e)}{2}  \\
        a_3  &= \frac{(a+d+f)}{2}  & \qquad
        a_4  &= \frac{(b+c+f) }{2} \\  
        b_1  &= \frac{(a+c + e +f)}{2}   & \qquad
        b_2  &= \frac{(b+d+e+f)}{2}  \\
        b_3  &= \frac{(a+b+c+d)}{2}  \\ 
        m &= min \lbrace a_1, a_2, a_3, a_4 \rbrace
        & \qquad M &= max \lbrace b_1, b_2, b_3 \rbrace 
\end{alignat*}

\begin{figure}[hbt] \epsfysize = 1.0 in
\centerline{\epsffile{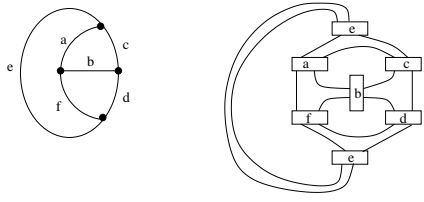}}
\caption{Tetrahedral net graph}
\label{fig:tetnet}
\end{figure}

The trivalent graph shown on the left in figure \ref{fig:tetnet} corresponds to the link diagram shown on the right.  The Jones polynomial of this virtual link diagram is
\begin{equation}
Tet \begin{bmatrix}  a & b & e \\ c & d & f \end{bmatrix}.
\end{equation} 

\begin{figure}[hbt] \epsfysize = 1.0 in
\centerline{\epsffile{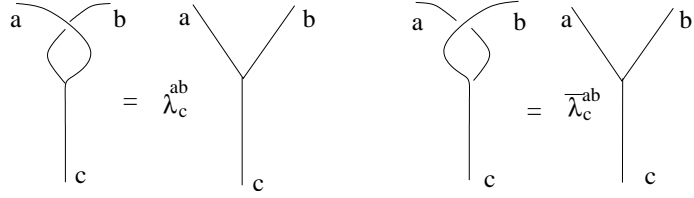}}
\caption{Graphical $ \lambda^{ab}_c $}
\label{fig:lamb}
\end{figure}

We may remove a classical crossing from a trivalent graph using the following formula.
We define the weight coefficient:
\begin{equation*}
  \lambda^{ab}_c =
  (-1)^{  \frac{a+b-c}{2} } A^{ \frac{a(a+2)+b(b+2)- c(c+2)}{2}}. 
\end{equation*}
The complex conjugate of  $  \lambda^{ab}_c $ is : 
\begin{equation}
 \bar{ \lambda}^{ab}_c =
  (-1)^{ \frac{a+b-c}{2} } A^{- \frac{a(a+2)+b(b+2)- c(c+2)}{2}}. 
\end{equation}
We remove a simple twist from the diagram as shown in figure \ref{fig:lamb}.
The equalities indicate that the bracket evaluations of the two diagrams are equivalent.

\begin{figure}[hbt] \epsfysize = 1.0 in
\centerline{\epsffile{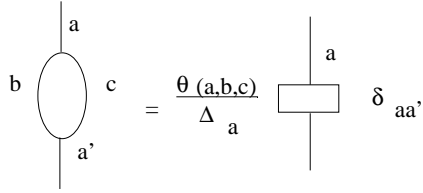}}
\caption{Bead removal formula}
\label{fig:beadrem}
\end{figure}

We may also simplify a trivalent graph using the formula shown in figure \ref{fig:beadrem}. We refer to this simplification as a bead removal.

\begin{figure}[hbt] \epsfysize = 1.0 in
\centerline{\epsffile{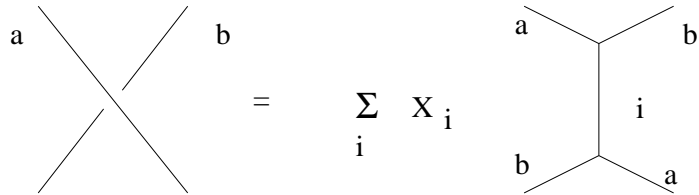}}
\caption{Crossing expansion formula}
\label{fig:crossx}
\end{figure}

The following formula equates a classical crossing with a weighted sum of trivalent subgraphs.
We exchange a crossing for a sum of trivalent subgraphs as shown in figure \ref{fig:crossx}, 

where
\begin{equation*}
X_i = \frac{\Delta_i \lambda^{ab}_i }{ \theta (a,b,i) }.
\end{equation*}

\subsection{Computations}

\begin{figure}[hbt] \epsfysize = 1.0 in
\centerline{\epsffile{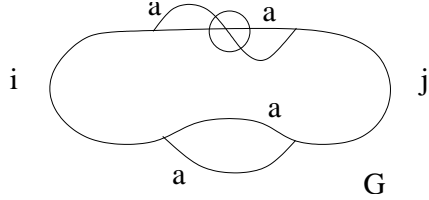}}
\caption{Graph obtained from $ K $ and $ \hat{K} $}
\label{fig:wrtform1}
\end{figure}

\begin{figure}[hbt] \epsfysize = 1.0 in
\centerline{\epsffile{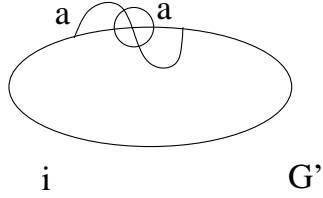}}
\caption{Reduced graph obtained from $ G $}
\label{fig:wrtform2}
\end{figure}

Let $ K $ and $ \hat{K} $ be the diagrams shown in figure \ref{fig:wrtex}.
We compute $ \langle K^{ \omega} \rangle $ and $ \langle \hat{K}^{ \omega} \rangle $.
Applying the definition of the WRT, 
\begin{align*}
 \langle K^{ \omega} \rangle &= \underset{a }{\sum} \Delta_a \langle K^a \rangle \\
 \langle \hat{K}^{ \omega} \rangle &= \underset{a }{\sum} \Delta_a \langle \hat{K}^a \rangle
 \end{align*}
Now, we use the formulas obtained from recoupling theory to
 transform the virtual knot 
diagrams $ K $ and $ \hat{K} $ into the trivalent graph $ G $ shown in figure \ref{fig:wrtform1}.
As a result: 
\begin{align*}
 \langle K^{ \omega } \rangle  &= \underset{i,j,a}{ \sum} \Delta_a \frac{ \Delta_i \Delta_j}{ \theta (a,a,i)
   \theta (a,a,j)} \lambda^{aa}_i \lambda^{aa}_j \langle  G  \rangle \\
   \langle \hat{K}^{ \omega } \rangle  &= \underset{i,j,a}{ \sum} \Delta_a  \frac{ \Delta_i \Delta_j}{ \theta (a,a,i)
   \theta (a,a,j)} \lambda^{aa}_i \lambda^{aa}_j  \lambda^{aa}_0 \langle  G  \rangle . 
\end{align*}
Note that $ (a,i,j) $ is not an admissible triple unless $ i=j $. Then: 
\begin{align*}   
   \langle K^{\omega} \rangle &= \underset{i,a}{ \sum} \Delta_a \frac{ \Delta_i^2 ( \lambda^{aa}_i )^2}{ \theta (a,a,i)^2 }   \langle G \rangle \\
  \langle \hat{K}^{\omega}  \rangle &= \underset{i,a}{ \sum} \Delta_a \frac{ \Delta_i^2 }{ \theta (a,a,i)^2}
   ( \lambda^{aa}_i )^2  \lambda^{aa}_0 \langle G \rangle
\end{align*}
Applying the $ \theta $-net identity, we obtain the graph, $ G' $, shown  
in  figure \ref{fig:wrtform2}. 

Now,
\begin{align*}
  \langle K^{ \omega} \rangle &= \underset{i,a}{ \sum} \Delta_a \frac{ \Delta_i }{ \theta (a,a,i) }
 ( \lambda^{aa}_i )^2 \langle  G' \rangle \\
 \langle \hat{K}^{\omega} \rangle &= \underset{i,a}{ \sum} \Delta_a \frac{ \Delta_i }{ \theta (a,a,i) }
 ( \lambda^{aa}_i )^2 \lambda^{aa}_0 \langle  G' \rangle\ 
\end{align*}

We determine for which values the triple, $ (a,a,i) $, is admissible. We then compute the Jones polynomial of the virtual link diagrams corresponding the the labeled graph, $ G' $. 

We indicate these values in table \ref{tab:JPR3}. If the triple $ (a,a,i ) $ is not admissible, we have entered a zero in the corresponding table entry.

\begin{table}
 \centering
 \begin{tabular}{||c|c|c||}
  \hline \textbf{ }  & $ a=0 $ & $ a=1 $  \\ \hline
  $ i=0  $   & $ 1$  & $ - A^2 - A^{-2} $  \\ \cline{1-3} \hline
  $  i=1 $  & $ 0 $  & $   0 $         \\  \hline
 \end{tabular}
 \label{tab:JPR3}
 \caption{Jones polynomial of the graph $ G' $, $ r=3$}
\end{table}

For $ r= 3 $:
\begin{gather*}
 \langle K^{\omega} \rangle =      \frac{ \Delta_0^2 ( \lambda^{00}_0)^2 \langle G' \rangle }{ \theta (0,0,0)}
   + \frac{ \Delta_0^2 ( \lambda^{11}_0)^2 \langle G' \rangle }{ \theta (1,1,0)}  
\end{gather*}
and
\begin{gather*}
\langle \hat{K}^{ \omega} \rangle =
\frac{ \Delta_0^2 ( \lambda^{00}_0)^3 \langle G' \rangle }{ \theta (0,0,0)}
   + \frac{ \Delta_0^2 ( \lambda^{11}_0)^3 \langle G' \rangle }{ \theta (1,1,0)}.   
 \end{gather*}
Now,
\begin{align*}
\langle K^{\omega} \rangle &= 0 \\
 \langle \hat{K}^{\omega} \rangle &= 1+ i
 \end{align*}
 
We fix $ r=4 $. We determine for what values $ (a,a,i) $ is an admissible triple. 
We then compute the Jones polynomial of the virtual link  diagram corresponding to the graph $ G_2 $.

\begin{table}
 \centering
 \begin{tabular}{||c|c|c|c||}
  \hline \textbf{ }  & $ a=0 $ & $ a=1 $ & $ a=2 $  \\ \hline
 $ i=0  $         &  $ 1 $ & $ -A^2 - A^{-2} $ & $ (-A^2 - A^{-2})^2 $      \\ \cline{1-4} \hline
 $ i=1 $           & $ 0 $   &  $      0   $    &  $    0   $               \\ \hline
 $ i=2  $              & $ 0 $  & $ -A^2-A^{-2} $  & $ 0 $  \\  \hline
 \end{tabular}
 \label{tab:JPR4}
 \caption{Jones polynomial of the graph $ G' $, $ r=4 $}
\end{table}

We compute $ \langle K^{ \omega } \rangle $ when $ r=4 $.
\begin{gather*}
\langle K^{\omega} \rangle =
 \frac{ \Delta_0^2 ( \lambda^{00}_0)^2 \langle G' \rangle }{ \theta (0,0,0)} 
  + \frac{ \Delta_0^2 ( \lambda^{11}_0)^2 \langle G' \rangle}{ \theta (1,1,0)} \\
  +  \frac{ \Delta_0^2 ( \lambda^{22}_0)^2 \langle G' \rangle }{ \theta (2,2,0)}   
  + \frac{ \Delta_2^2 ( \lambda^{11}_2)^2 \langle G' \rangle}{ \theta (1,1,2)}  
\end{gather*}
Similarly, we compute the WRT of $ \hat{K} $ when $ r=4 $.
\begin{gather*} 
\langle \hat{K}^{\omega} \rangle =
\frac{ \Delta_0^2 ( \lambda^{00}_0)^3 \langle G' \rangle }{ \theta (0,0,0)} 
  + \frac{ \Delta_0^2 ( \lambda^{11}_0)^3 \langle G' \rangle}{ \theta (1,1,0)} \\
  +  \frac{ \Delta_0^2 ( \lambda^{22}_0)^3 \langle G' \rangle }{ \theta (2,2,0)}   
  + \frac{ \Delta_2^2 ( \lambda^{11}_2)^2 \lambda^{11}_0 \langle G' \rangle}{ \theta (1,1,2)}  
\end{gather*}
We determine that
\begin{align*}
\langle K \rangle &= 1.29289 + 1.70711 i \\
\langle \hat{K} \rangle &= 1.23044 + 0.92388 i
\end{align*}
The virtual knot diagram $ K $ has 
 $ b_{+} (K) =0 $ and $ b_{-} (K) =0 $ implying $ n(K) =-1$. 
Now, we normalize the WRT and determine that:
\begin{align*}
Z_K (3) &= 0 \\
Z_{\hat{K}} (3) &= 0.707107 i
\end{align*}
We determine that the normalized WRT when $ r=4 $.
\begin{align*}
Z_K (4) &= -0.517982 + 0.135299 i \\
Z_{\hat{K}}(4) &= -0.331106 + 0.195807 i
\end{align*}

\section{Recoupling Theory and Virtual Knot Diagrams} \label{recoup}

Recall the recoupling theorem \cite{tl}.
\begin{thm}[Recoupling Theorem] 
A single labeled trivalent graph is equivalent to a sum of labeled trivalent graphs with coefficients in $ \mathbb{C} $ as illustrated in figure \ref{fig:graphx}.
\begin{figure}[hbt] \epsfysize = 1.0 in
\centerline{\epsffile{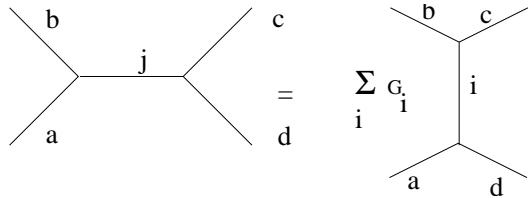}}
\caption{Graph exchange (exchanging a H-shaped graph for an I-shaped graph)}
\label{fig:graphx}
\end{figure}
Note that
\begin{equation}
G_i = \frac{ Tet \begin{bmatrix} a & b & i \\
												c & d & j \end{bmatrix} \Delta_i }{ \theta (a,d,i) \theta (b,c,j) }
\end{equation}
\end{thm}
\textbf{Proof:} \cite{tl}. \qed

\begin{figure}[hbt] \epsfysize = 1.0 in
\centerline{\epsffile{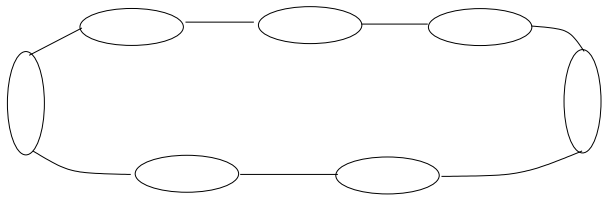}}
\caption{Example: String of Beads Graph}
\label{fig:trivalent}
\end{figure}

\begin{rem} One of the formulas given in section \ref{wrti} is a special case of the recoupling theorem.
In particular, if $ j=0 $, we obtain the fusion formula from section \ref{wrti}. \end{rem}
We now prove a theorem that demonstrates a sum of $ H $-shaped trivalent graphs can be rewritten as a sum of $ I $ shaped trivalent subgraphs. This procedure gives the appearance of sliding edges past the vertices of the trivalent graph. Using the recoupling theorem, we 
may transform any trivalent graph without virtual crossings into a sum of weighted \emph{string of beads} graphs as illustrated in figure \ref{fig:trivalent}.

\begin{lem}[Trivalent Graphs]
A classical knot diagram may be transformed into a string of beads  
graph by replacing each classical crossing with a trivalent subgraph and applying the recoupling theorem.
\end{lem}
\begin{rem} This result follows from the fact that thickening an I-shaped and a H-shaped graph produces two homeomorphic surfaces. The surface produced by thickening any planar graph is a n-punctured sphere. The punctured sphere deformation retracts to a string of beads graph after homeomorphism. We will generalize this approach to show that virtual link diagrams correspond to a specific type of trivalent graphs. This lemma will be proved as a sub case of Lemma \ref{virtgraph}. Refer to \cite{tl}, p. 115.
\end{rem}

The recoupling formula allows us to translate any planar trivalent graph (no virtual crossings) into a string of beads graph. This fact restricts the types of subgraphs that occur (in graphs obtained from classical link diagrams) to a small number of possibilities. As a result, we need only compute formulas that determine the WRT of the classical tangle diagrams that correspond these labeled subgraphs.
In particular, it is sufficient to compute the coefficients that correspond to the insertion and deletion of 
$ \theta $-net subgraphs. 

\begin{figure}[hbt] \epsfysize = 0.75 in
\centerline{\epsffile{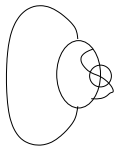}}
\caption{Unreduceable super bead}
\label{fig:twistgraph}
\end{figure}

The reduction to a simple string of beads does not hold for virtual link diagrams. We can remove all classical crossing from the diagram by exchanging the classical crossings for H or I shaped trivalent subgraphs, and obtain a trivalent graph with only virtual crossings. However, we can not remove all virtual crossings by applying the recoupling theorem. 
For example, the graph shown in figure \ref{fig:twistgraph} can not be reduced to a trivalent graph without virtual crossings. 

\begin{figure}[hbt] \epsfysize = 0.75 in
\centerline{\epsffile{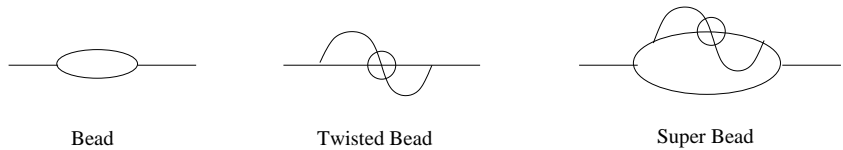}}
\caption{Beads types formed by virtual knot diagrams}
\label{fig:bead}
\end{figure}

The reduced form of a virtual knot diagram is similar to a string of beads graph but includes three different types of beads. The reduced form of the virtual knot diagram is a string of  \emph{super beads}, \emph{twisted beads} and \emph{beads}. We illustrate the possible three  bead types from the virtual case in figure \ref{fig:bead}.
 
\begin{figure}[hbt] \epsfysize = 1.0 in
\centerline{\epsffile{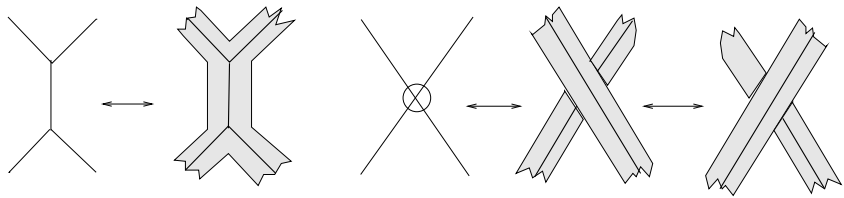}}
\caption{Thickened graph}
\label{fig:thick}
\end{figure}

We prove that the WRT of a virtual link diagram is equal to a weighted sum of the Jones polynomial of
virtual link diagrams. The virtual link diagrams are indicated by labeled 
virtual string of beads graphs. These string of bead graphs are obtained from the virtual link diagram by applying the recoupling theorem and crossing exchange.
We first prove a lemma that demonstrates any trivalent graph with virtual crossings can be transformed into 
a string of the three bead types by sliding edges.
\begin{lem}[Graphical Lemma]\label{virtgraph} Any trivalent
 graph with virtual crossings may be 
reduced to a trivalent string of beads graph with three possible types
of beads: the standard bead, twisted bead, and super bead. This graph will contain at most one twisted bead. 
\end{lem}
\textbf{Proof:} Thicken the graph so that the graph is
  the core of a compact oriented 
surface, $ T(G) $. An example of a graph and its corresponding thickened surface is
shown in figure \ref{fig:thick}. Note that we may apply homeomorphisms these surfaces and the surfaces do not depend on the choice of embedding in 3 dimensional space.
The two surfaces on the right hand side of figure \ref{fig:thick} are equivalent.

Note that if the graph, $ G $  has $ 2n $ vertices, then $ G $ has $ 3n $ edges
since the graph is trivalent. The Euler characteristic of the surface, $ T(G) $
 is $ -n $ (since G is a deformation retract of T(G)). Without loss of generality, we assume that the surface
 has $ k $ boundary components.

Recall the following theorem from \cite{hirsch}.
\begin{thm}[Surface Classification] \label{hir}
 Let $ M $ be a connected, compact, orientable surface of
Euler characteristic $ \chi $. Suppose $ \partial M $ has $ k \geq 0 $ 
boundary components. Then $ \chi + k $ is even. Let
 $ p = 1 - \frac{ (\chi +k)}{2} $ then M is diffeomorphic to the surface 
obtained from an orientable surface of genus $ p $ by the removal of the interiors of
k disjoint disks. 
\end{thm}

\begin{figure}[hbt] \epsfysize = 1.0 in
\centerline{\epsffile{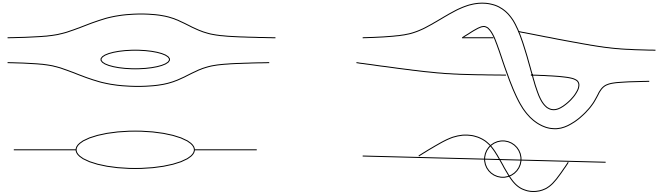}}
\caption{Surfaces corresponding to the standard bead and the twisted bead}
\label{fig:surfgraph2}
\end{figure}

\begin{figure}[hbt] \epsfysize = 1.0 in
\centerline{\epsffile{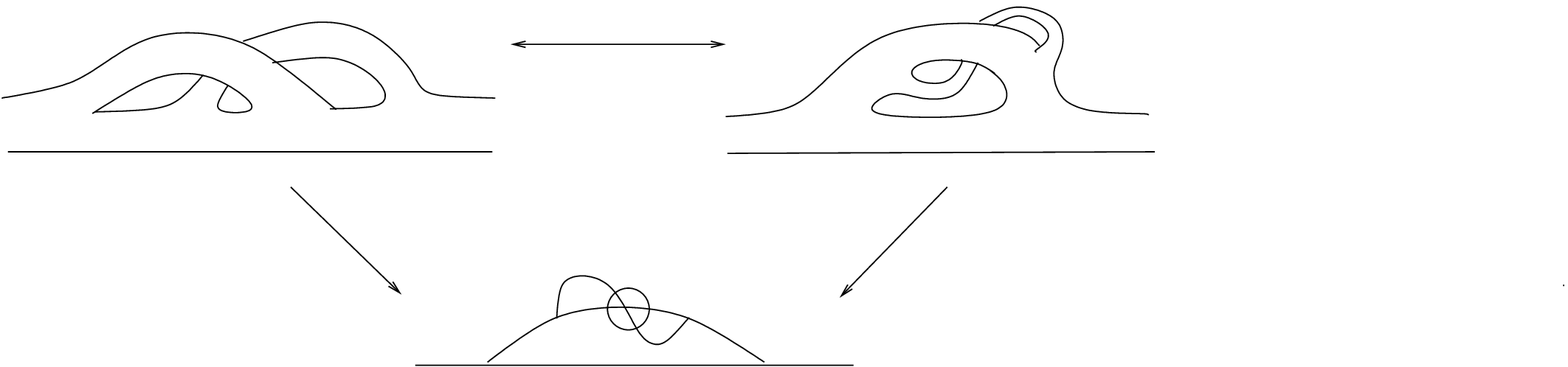}}
\caption{Surfaces corresponding to the super bead}
\label{fig:surfgraph3}
\end{figure}

Applying Theorem \ref{hir}, $ \chi (S) + k = -n + k $ and $ -n + k $ is 
even. Let $ p = 1 - \frac{ -n + k}{2} $. We observe that $ T(G) $ is the connected sum 
of $ p $ tori with $ k $ punctures. 
Each torus component contributes a handle pair. Each handle pair corresponds to a either a twisted bead or part of a super bead depending on the number of punctures. Each punctured disk forms a either plain bead or part of a super bead as shown in figures \ref{fig:surfgraph2} and \ref{fig:surfgraph3}.

\begin{figure}[hbt] \epsfysize = 1.75 in
\centerline{\epsffile{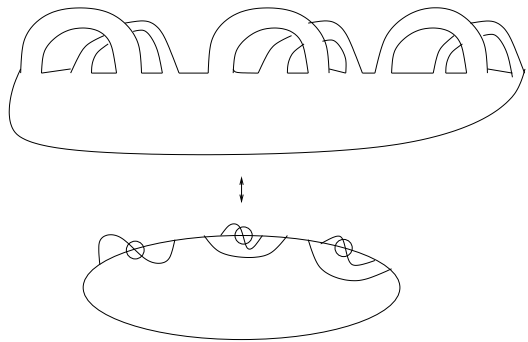}}
\caption{A surface  and its  corresponding graph}
\label{fig:surfgraph1}
\end{figure}

For example, a surface with one handle pair and one puncture corresponds 
to a single twisted bead on a loop (the graph shown in figure \ref{fig:twistgraph}). 
The connected sum of three tori with one puncture and its corresponding graph are shown in figure \ref{fig:surfgraph1}.

After isotopy,
any surface that is the connected sum of $ p $ tori
with $ k $ punctures as may be regarded as a series of $ p $ handle pairs attached to a disk with $ k-1 $ 
punctures.  
As a result, 
 these surfaces deformation retract  to a string of beads type trivalent  graph with standard beads, twisted
beads, and a single super bead.  \qed

\begin{thm}[Virtual Recoupling Theorem] 
Let $ K $ be a labeled virtual link diagram. Using the recoupling theorem and crossing expansion, we may transform the virtual link diagram into a sum of labeled
trivalent, string of beads graphs with coefficients in $ \mathbb{C} $.
The WRT of the labeled virtual link diagram is equal to the WRT of the sum of 
weighted trivalent graphs.
\end{thm} 
\textbf{Proof:} From the previous theorem, it is clear that 
we can transform a virtual link diagram into a 
string of beads graph using the 
recoupling theorem  and crossing expansion without coefficients.
After determining this sequence of transformations, apply the 
appropriate coefficients and labels indicated by the recoupling theorem, the crossing expansion formula, 
and the initial labeling of the link. \qed

\subsection{Obstruction to Weight Formulas} \label{brobs}
Recall that the $ n^{th} $ Jones-Wenzl idempotent has the property that $ T_n U_i = 0 $ for any elementary n-tangle $ U_i $. 
We cannot easily extend the evaluation formulas given earlier in this paper because 
 the Jones-Wenzl projector does not necessarily have the property that $ T_n P_n = 0 $ when $ P_n $ is a n-tangle with only virtual crossings.
 In 
 the classical case, we relied on this fact to simplify the recursive formulas. The recursion formulas needed to determine the value of the Jones polynomial of $ T_n P_n $ depend on  
the position of the virtual crossings. The evaluation of these formulas 
requires that we know the Jones polynomial
of the closure of $ T_{n-1} P_{n-1} $ for all $ (n-1)$-tangles $ P_{n-1} $ containing only virtual crossings. 
We prove two theorems that indicate the complexity involved in computing the Jones polynomial of $ T_n P_n $

\begin{thm}\label{noidem} The $ n^{th} $ J-W projector, $ T_n $, does not annihilate all n-tangles $ P_n $ with only virtual crossings. 
\end{thm}
To prove Theorem \ref{noidem} we need only demonstrate that
$ T(T_n E) \neq 0 $ where $ E $ is a n-tangle with a single virtual crossing.
We postpone this proof until after the proof of Lemma \ref{lemnotzero}.

\begin{figure}[hbt] \epsfysize = 2.0 in
\centerline{\epsffile{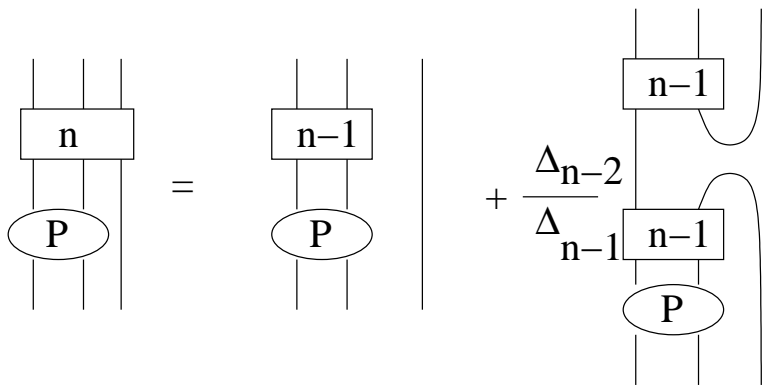}}
\caption{Diagrammatic proof of Proposition \ref{pp1}}
\label{fig:perm1proof}
\end{figure}

We introduce the following propositions preliminary to 
the proof of Theorem \ref{noidem}.
The formulas in these propositions determine the number of the Jones polynomial of the closure $ T_n P_n $.
\begin{prop} \label{pp1}
\label{perm1} If $ P $
is a (n-1)-tangle with only virtual crossings, we construct a n-tangle, $ P_n^* $, by attaching an additional strand on the right. 
Then the Jones polynomial of the closure of $ T_n P^* $ is determined by the following formula.
$$
T_{n} P = T_{n-1} P - \frac{ \Delta_{n-2} }{ \Delta_{n-1} } T_{n-1} P
$$
\end{prop}
\textbf{Proof: } Apply the recursion relation shown in figure \ref{fig:perm1proof}. \qed

\begin{prop} \label{pp2}
Let $ P $ be a n-tangle with only virtual crossings. We compute 
$ T_n P $ recursively, and the depth of the recursion 
equation depends on the initial position of the rightmost strand in the n-tangle $ P $. 
\end{prop}

\begin{figure}[hbt] \epsfysize = 3 in
\centerline{\epsffile{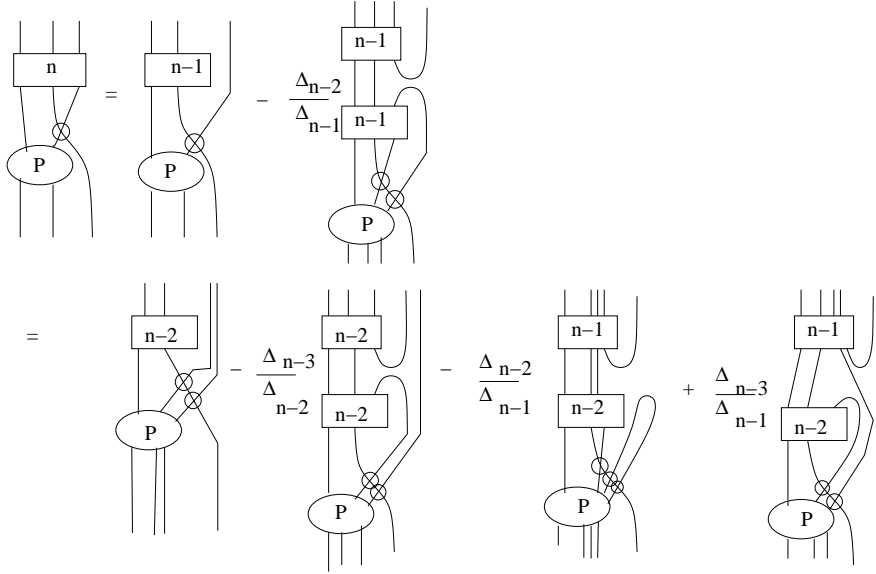}}
\caption{Diagrammatic proof of Proposition \ref{pp2}}
\label{fig:perm2proof}
\end{figure}

\textbf{Proof:}
Let $ P$ be an n-tangle with only virtual crossings. We compute $ T_n P $ recursively as shown in figures \ref{fig:perm1proof} and \ref{fig:perm2proof}.
We expand the diagram as in indicated in figure \ref{fig:perm2proof}. 
The complexity of the recursion depends on the virtual crossings in $ P $. \qed
\begin{rem} The question of handling this recursion in a systematic way is an interesting research question.
\end{rem}
\begin{lem} \label{lemnotzero} 
Let $ E $ denote the n-tangle with a single virtual crossing involving only the $ 1^{st} $ and the $ 2^{nd} $ strands.  Then:
 $$
 d \langle cl(T_n E) \rangle =(d + \frac{ \Delta_{n-2}}{ \Delta_{n-1} }) \dots (d +
 \frac{ \Delta_{1}}{ \Delta_{2} }) (d-1)
$$ 
\end{lem}
\textbf{Proof:} 
We apply formula from Proposition \ref{perm1} to $ T_n E $.
We calculate $ \langle cl( T _2 E ) \rangle  = 1-1/d $. 
 \qed

We now prove Theorem \ref{noidem}.

\textbf{Proof of Theorem \ref{noidem}: }
If $ T_n $ is an idempotent, then $ T_n E = 0 $. However,  $ \langle cl(T_n E) \rangle \neq 0 $ 
indicating that $ T_n E $ is not the empty diagram. \qed

The algebra generated by  extending n-tangle diagrams by including virtual crossings is a diagrammatic version of the Brauer algebra. Note that each classical n-tangle diagram is expanded (via the bracket) into an element of the Temperly-Lieb algebra. It is this extension of the Temperly-Lieb algebra by virtual elements that is isomorphic to the Brauer alegebra. The Brauer algebra is described in \cite{benk}. We would like to understand the structure of the Jones-Wenzl projectors in the context of the Brauer algebra. Results in this area would help in understanding the extension of the WRT invariant that we have discussed in this paper.

\noindent {\bf Acknowledgement.} This effort was sponsored in part by the
Office of the Dean at the United States Military Academy.
 The U.S. Government is authorized to reproduce and distribute
reprints
for Government purposes notwithstanding any copyright annotations thereon.
The
views and conclusions contained herein are those of the authors and should
not be
interpreted as necessarily representing the official policies or
endorsements,
either expressed or implied, of the United States Military Academy or the U.S. Government. (Copyright
2004.)
Much of this effort was sponsored by the Defense Advanced Research Projects Agency (DARPA) and Air Force Research Laboatory, Air Force Materiel Command, USAF, under agreement F30602-01-2-05022. The US Government is authorized to reproduce and distribute reprints for Government purposes notwithstanding any copyright annotations thereon. The views and conclusions contained herein are those of the authors and should not be either expressed or implied, of the Defense Advanced Research Projects Agency, the Air Force Research Laboratory, or the U.S. Government. (Copyright 2004.) It gives the second author great pleasure to acknowledge support from NSF Grant DMS-0245588, and to thank the Stanford Linear Accelerator Theory Group, the University of Waterloo and the Perimeter Institute for hospitality during the preparation of this research.


\begin{thebibliography}{10}

\baselineskip=12pt 
\parskip=2pt plus 1pt 



\bibitem{benk}
 \newblock{Georgia Benkart, }
 \newblock{Commuting actions---a tale of two groups, Contemporary Mathematics, }  
\newblock{American Mathematical Society, Vol. 194, p. 1-46, 1996, }
\newblock{Lie algebras and their representations (Seoul, 1995)} 

\bibitem{cks}
\newblock{J. Scott Carter, Seiichi Kamada, and Masahico Saito, }
\newblock{ Stable equivalence of knots on surfaces and virtual
  knot cobordisms, }
\newblock{ Knots 2000 Korea, Vol. 1 (Yongpyong), }
\newblock{ J. Knot Theory Ramifications 11 (2002), no. 3, 311--32}


\bibitem{kd1}
 \newblock{Heather A. Dye and Louis H. Kauffman, }
\newblock{Minimal surface representations of virtual knots and links, }
 \newblock{www.arxiv.org, GT/0401035}

\bibitem{GPV}
\newblock{Mikhail Gussarov, Michael Polyak, and Oleg Viro}
\newblock{Finite Type Invariants of Classical and Virtual Knots, }
\newblock{Toplogy 39 (2000), no. 5 p. 1045-1068}

\bibitem{hirsch}
\newblock{Morris W. Hirsch, }
\newblock{Differential Topology, }
\newblock{Springer-Verlag, Graduate Texts in Mathematics, 1997}

\bibitem{kamada}
\newblock{Naoko Kamada and Seiichi Kamada, }
\newblock{Abstract link diagrams and virtual Knots, }
\newblock{Journal of Knot Theory and it's Ramifications, Vol. 9 No. 1, p. 93-109, }   
\newblock{World Sci. Publishing, 2000}

\bibitem{Kdetect}
\newblock{Louis H. Kauffman, }
\newblock{Detecting Virtual Knots, }
\newblock{Atti del Seminario Matematico e Fisico
         dell'Universite di Modena, Vol. 49, suppl., p. 241-282,}
\newblock{Univ. Modena, 2001}

\bibitem{knotphys}
\newblock{Louis H. Kauffman, }
\newblock{Knots and Physics, }
\newblock{Series on Knots and Everything, Vol. 1, World Scientific, 1991, 1994, 2001}

\bibitem{kvirt}
\newblock{Louis H. Kauffman, }
\newblock{Virtual Knot Theory, } 
\newblock{European Journal of Combinatorics, Vol. 20, No. 7, p. 663-690,}
\newblock{Academic Press, 1999}



\bibitem{tl}
\newblock{Louis H. Kauffman and Sostenes L. Lins, }
\newblock{Temperly-Lieb Recoupling Theory and Invariants of 3-Manifolds, }
\newblock{Annals of Mathematics Studies, Princeton University Press, 1994 }




\bibitem{kup}
\newblock{Greg Kuperberg, }
\newblock{What is a Virtual Link?,}
\newblock{Algebr. Geom. Topol. 3 (2003), p. 587-591 (electronic)}
\newblock{www.arxiv.org, GT/0208039}


\bibitem{vom}
\newblock{Vassily O. Manturov, }
\newblock{Kauffman-Like Polynomial and Curves in 2-Surfaces}
\newblock{Journal of Knot Theory and its Ramifications, Vol. 12 (2003), no. 8, p. 1145-1153}

\bibitem{purp}
\newblock{Viktor V. Prasolov and Alexei B. Sossinsky, }
\newblock{ Knots, Links, Braids and 3-Manifolds;
 An Introduction to New Invariants in Low-Dimensional Topology, }
\newblock{American Mathematical Society, Translations of Mathematical
 Monographs, 1996}

 \bibitem{rt1}
 \newblock{Nikolai Yu Reshetikhin and Vladimir G. Turaev, }
 \newblock{Invariants of 3-manifolds via link polynomials and quantum groups, }
 \newblock{Inventiones Mathematicae, Vol. 103, No. 3, p. 547-597 },
 \newblock{Springer, 1991}


 \bibitem{rt2}
 \newblock{Nikolai Yu Reshetikhin and Vladimir G. Turaev, }
 \newblock{Ribbon Graphs and their invariants derived from quantum groups, } 
 \newblock{Communications in Mathematical Physics, Vol. 127, No. 1, p. 1-26, }
 \newblock{Springer, 1990},






\bibitem{witten}
\newblock{Edward Witten,}
\newblock{Quantum field theory and the Jones polynomial, }
\newblock{Braid group, knot theory and statistical mechanics, II, p. 361-451, }
\newblock{ Adv. Ser. Math. Phys., 17, World Sci. Publishing 1994}
\end{thebibliography}
\end{document}